\documentclass{conm-p-l}
\usepackage{amsmath,amsthm,amssymb, amstext, amsopn, amsxtra}

\theoremstyle{definition}

\theoremstyle{remark}

\numberwithin{equation}{section}
\newcommand{\nt}{\noindent}
\newcommand{\bs}{\bigskip}

\newcommand{\mk}{\medskip}
\newcommand{\sk}{\smallskip}
\newcommand{\D}{\Delta}
\newcommand{\G}{\Gamma}

\begin{document}

\title{Wedderburn Polynomials over Division Rings}

\author{T. Y. Lam}
\address{Department of Mathematics, University of California, 
Berkeley, CA 94720}
\email{lam@math.berkeley.edu}
\thanks{Some of the results in this note were presented by the first
author in a lecture (with the same title) at the National Conference on 
Algebra VII, at Beijing Normal University in October, 1999. We owe a 
debt of gratitude to Professor Zhang Yingbo, who kindly provided a Chinese
translation of an English version of this article.}
\author{Andr\'e Leroy}
\address{Department of Mathematics, Universit\'e d'Artois, 62307 Lens Cedex,
France}
\email{leroy@euler.univ-artois.fr}

\subjclass{Primary 16D40, 16E20, 16L30; Secondary 16D70, 16E10, 16G30}
\date{January 3, 2000.}

\begin{abstract}
A Wedderburn polynomial over a division ring $K$ is a minimal polynomial
of an algebraic subset of $K$.  Special cases of such polynomials include,
for instance, the minimal polynomials (over the center $\,F=Z(K)$) of elements 
of $\,K\,$ that are algebraic over $\,F$.  In this note, we give a  
survey on some of our ongoing work on the structure theory of Wedderburn 
polynomials.  Throughout the note, we work in the general setting of an 
Ore skew polynomial ring $\,K[t,S,D]$.
\end{abstract}

\maketitle

\begin{center}
{\bf \S1. \  Introduction}
\end{center}

\medskip
The main object of this paper is the study of a class of polynomials over a 
division ring $K$, which we call {\it Wedderburn polynomials\/} (or simply 
W-polynomials).  Roughly speaking, a W-polynomial over $K$ is one which has 
``enough zeros'' in $K$.  (For a more precise definition, see (3.1).)  In 
the case when $K$ is a field, W-polynomials are simply those of the form 
$\,(t-c_1)\cdots (t-c_n)$, where $\,c_1,\dots,c_n\,$ are {\it distinct\/} 
elements of $\,K$. In the general case of a division ring $\,K$, a 
W-polynomial still has the form $(t-c_1)\cdots (t-c_n)$, although the 
$\,c_i\,$'s need no longer be distinct.  And even if the $\,c_i\,$'s are 
distinct, $\,(t-c_1)\cdots (t-c_n)\,$ need not be a W-polynomial.  The 
recognition of a W-polynomial, for instance, is already a very interesting 
problem over a division ring $\,K$.

\bs

The early work of Wedderburn [We] (ca.~1921) showed that, if $\,a\in K\,$ 
is an algebraic element over the center $\,F\,$ of $\,K$, then the minimal 
polynomial of $\,a\,$ over $F$ (in the usual sense) is a W-polynomial in 
$\,K[t]$ (and in particular splits completely over $\,K$).  This classical 
result of Wedderburn has led to much research on $\,K[t]$, and has found 
important applications to the study of subgroups of the multiplicative group 
$\,K^*$, central simple algebras of low degrees and crossed product algebras,
PI-theory, Vandermonde matrices, Hilbert 90 Theorems, and the theory of 
ordered division rings, etc. For some literature along these lines, see 
[Al], [HR], [Ja$_3$], [La], [LL$_1$]-[LL$_3$], [Ro$_1$]-[Ro$_3$], [Se], 
[RS], ...\,, etc.

\bs
Our definition of W-polynomials was directly inspired by the afore-mentioned 
work of Wedderburn, although our W-polynomials will have coefficients in 
$\,K$, instead of in $\,F$.  These W-polynomials have different properties 
and characterizations, and seem to be quite basic in working with the 
polynomial theory over $\,K$. We have recently initiated a systematic
investigation of the theory of W-polynomials.  In this note, we'll give a 
synopsis of some of the results obtained so far.  These include, for
instance:

\begin{itemize}
\item A theorem on the left-right symmetry of W-polynomials;

\sk
\item Several criteria for the factors and products of W-polynomials to be 
W-polynomials;

\sk
\item Some ``Rank Theorems'' describing the behavior of the rank 
of algebraic sets with respect to union, intersection, and certain other
transformations;

\sk
\item Applications of quadratic W-polynomials to the study of certain 
``metro-equations'' in division rings; and

\sk
\item A study of the relationship between the set of W-polynomials and 
the set of algebraic subsets of $\,K$, from the viewpoint of complete
modular lattices.
\end{itemize}

\mk
Following Ore [Or], we work in the setting of skew polynomials (rather than 
just ordinary polynomials) over the division ring $\,K$.  This added degree 
of generality is definitely worthwhile considering that skew polynomials 
have become increasingly important with the growing interests in quantized 
structures and noncommutative geometry.  The basic mechanism of skew 
polynomials is recalled in \S2, where we also set up the terminology and 
general framework for the paper. Wedderburn polynomials are introduced in
\S3, and the synopsis of the main results on W-polynomials is given in 
\S\S\,4-8.

\bs
Details of the research reported in this note (including the proofs of
the results indicated above) will appear in the authors' forthcoming papers 
[LL$_4$] and [LL$_5$].

\sk

\bs
\begin{center}
{\bf \S2. \  Skew Polynomial Rings and Algebraic Sets in $\,K$}
\end{center}

\bs
To work with skew polynomials, we start with a triple $(K,S,D)$, where $K$ 
is a division ring, $S$ is a ring endomorphism of $K$, and $D$ is an 
$S$-derivation on $K$.  (The latter means that $D$ is an additive endomorphism
of $K$ such that $D(ab)=S(a)D(b)+D(a)b, \; \forall a,b \in K$.)  In this 
general setting, we can form $\,K[t,S,D]$, the {\it Ore skew polynomial 
ring\/} consisting of polynomials of the form 
$\,\sum b_{i}t^{i}$ \ ($b_{i} \in K$) which are added in the usual way 
and multiplied according to the rule 
$$\,tb=S(b)t+D(b) \;\;\;(\,\forall \, b \in K). \leqno (2.1) $$
In case $\,(S,D)=(I,0)\,$ (we'll refer to this as the ``classical case''), 
$\,K[t,S,D]\,$ boils down to the usual polynomial ring $\,K[t]\,$ with a 
{\it central\/} indeterminate $\,t$.  {\it Throughout this paper, we'll write 
$\,R:=K[t,S,D]$.}  It is easy to check that $\,R\,$ admits an euclidean 
algorithm for right division, so $\,R\,$ is a principal left ideal domain.

\bs
In working with $\,R$, it is important to be able to ``evaluate''
a polynomial $\,f(t)= \sum b_{i}t^{i}\,$ at any scalar $\,a\in K$, that is,
to define $\,f(a)$. Following our earlier work [LL$_1$], we take $\,f(a)\,$
to be $\,\sum b_{i}N_{i}(a)$, where the ``$i\,$th power function'' $\,N_i\,$
is defined inductively by
$$  N_{0}(a)=1,  \; \; \; \mbox{and} \; \; \; 
N_{i}(a)=S(N_{i-1}(a))a+D(N_{i-1}(a)), \; \; \;\forall \,a \in K.\leqno (2.2)$$
That this gives the ``right'' definition of $\,f(a)\,$ is seen from the
validity of the {\it Remainder Theorem\/} [LL$_1$:~(2.4)]: there is a unique 
$\,q\in R$ such that 
$$\,f(t)=q(t)(t-a)+f(a).  $$
From this, it follows 
immediately that $\,f(a)=0\,$ iff $\,t-a\,$ is a right factor of $\,f(t)$.  

\bs
Another remarkable fact about evaluating skew polynomials at scalars is the
``Product Formula'' [LL$_1$:~(2.7)] for evaluating $\,f=gh\,$ at any 
$\,a\in K$:
$$ (gh)(a)= \left\{ \begin{array}{ll}
0                 & \; \;\; \mbox{if} \; \;\;h(a)=0,  \\
g(a^{h(a)})h(a)   & \; \;\; \mbox{if} \; \;\;h(a)\neq 0.
\end{array}
\right.     \leqno (2.3) $$
Here, for any $\,c\in K^{*}$, $\,a^{c}\,$ denotes $\,S(c)ac^{-1}+D(c)c^{-1}$,
which is called the $\,(S,D)$-{\it conjugate\/} of $\,a\,$ (by $\,c$).  With 
this general conjugation notation, it is easy to verify by a direct 
calculation that
$$ (a^c)^d=a^{dc} \;\;\; \mbox{for\;any} \;\,c,\,d\in K^*. \leqno (2.4) $$
From this, it follows readily that $\,(S,D)$-conjugacy is an equivalence
relation.  In the following, we shall write 
$$\,\Delta^{S,D}(a):= \{\,a^{c}: \;c\in K^{*}\,\};  \leqno (2.5)$$
this is called the {\it $\,(S,D)$-conjugacy class\/} of $\,a$. All such 
classes form a partition of $\,K$.  For instance, $\,\Delta^{S,D}(0)\,$ is the
set of all {\it logarithmic derivatives} $\,\{D(c)c^{-1}:\; c\in K^{*}\}$.  
And, in the classical case, $\,\Delta^{I,0}(a)\,$ is just the ``usual'' 
conjugacy class 
$$\,\D(a)=\{cac^{-1}:\;c\in K^*\}.  $$  

\mk
Next, we introduce a basic notation for this paper.  For $\,f\in R$, let 
$$\,V(f):=\{a\in K:\; f(a)=0\}. \leqno (2.6)$$
We'll say that a set $\,\D\subseteq K\,$ is {\it algebraic\/} (or, more 
precisely, $\,(S,D)$-algebraic) if $\,\D\subseteq V(f)\,$ for some
nonzero $\,f\in R$.  In this case, the set of polynomials vanishing 
on $\,\D\,$ forms a nonzero left ideal in $\,R$. The monic generator of this
left ideal is called the {\it minimal polynomial\/} of $\,\D$; we denote 
it by $\,f_{\D}$.  The degree of $\,f_{\D}\,$ is called the {\it rank\/} 
of the algebraic set $\,\D\,$; we denote it by \,rk$\,(\D)$.   According to 
the Remainder Theorem, $\,f_{\D}\,$ is just the (monic) ``llcm'' (least left 
common multiple) of the linear polynomials $\,\{t-a:\;a\in \D\}$.  As in 
[La:~Lemma 5], it is easy to see that $\,f_{\D}\,$ has always the form 
$\,(t-c_1)\cdots (t-c_n)\,$ where each $\,c_i\,$ is $\,(S,D)$-conjugate to 
some element of $\,\D$. 

\bs
Of course, all of the above was inspired in part by classical algebraic
geometry.  Going a little further, we get a theory of {\it polynomial 
dependence\/} (or P-dependence for short) for the elements of $\,K$.  By 
definition, an element $\,d\in K\,$ is P-dependent on an algebraic set 
$\,\D\,$ if $\,g(d)=0\,$ for every $\,g\in R\,$ vanishing on $\,\D$.  We see 
easily that the set of elements P-dependent on $\,\D\,$ is precisely 
$\,V(f_{\D})$, which we shall henceforth call the ``P-closure'' of $\,\D\,$ 
and denote by $\,\overline{\D}$. As in [La], we can also define 
{\it P-independence\/} and the notion of a {\it P-basis\/} for an algebraic 
set $\,\D\,$ in a natural manner.  The cardinality of a P-basis for $\,\D\,$ 
is just rk$\,(\D)$.  If $\,\{a_1,\dots,a_n\}\,$ is a P-basis of $\,\D\,$
(where $\,n=\mbox{deg}(f_{\D})$), then $\,f_{\D}\,$ is in fact the llcm of 
the linear polynomials $\,\{t-a_i:\,1\leq i\leq n\}$.  We refer the reader 
to [La] for the rudiments of the theory of P-dependence. (Although this 
theory was developed in [La] in the case $\,D=0$, it holds word-for-word 
also in the $\,(S,D)$-case.)

\bs

\sk
\begin{center}
{\bf \S3. \  The Class $\,{\mathcal W}\,$ of Wedderburn Polynomials}
\end{center}

\bs
We now come to the principal object of the paper.

\bs\nt
{\bf Definition 3.1.}  A polynomial $\,f\in R\,$ is said to be a 
{\it Wedderburn polynomial\/} (or simply a \,W-polynomial) if it is the
minimal polynomial of some algebraic set in $\,K$, or, equivalently,
if $\,f=f_{V(f)}$. The set of W-polynomials in $\,R\,$ will be denoted
throughout this paper by the symbol $\,{\mathcal W}$.

\bs
Some easy characterizations of W-polynomials are given in the following
result.

\bs\nt
{\bf Proposition 3.2.} {\it For a monic polynomial $\,f\in R\,$ of degree 
$\,n$, the following are equivalent\/}:

\sk
(1) {\it $f$ is a \,W-polynomial\/};  

\sk
(2) rk$\,(V(f))=n\,$ (``$f\,$ {\it has enough zeros\/}'');  

\sk
\indent (3) {\it For any $\,p\in R$, $\,V(f)\subseteq V(p) \Longrightarrow 
p\in R\cdot f$.}

\bs
To motivate the study of W-polynomials, we mention the following 
basic examples.

\bs\nt
{\bf Examples 3.3.} 

\mk
(1) The constant polynomial $\,1\in R\,$ is a W-polynomial, as it is
the minimal polynomial of the empty (algebraic) set $\,\emptyset$.

\sk
(2) For any $\,a\in K$, $\,t-a\,$ is a W-polynomial, as it is the minimal
polynomial of the singleton (algebraic) set $\,\{a\}$.

\sk
(3) An easy argument shows that a monic quadratic polynomial $\,f\in R\,$
is a W-polynomial iff $\,|V(f)|\geq 2$.  Alternatively, the quadratic
W-polynomials are those of the form $\,\bigl(t-b^{b-a}\bigr)(t-a)$,
where $\,a\neq b\,$ in $\,K$. For instance, over the quaternions,
with $\,(S,D)=(I,0)$, $\,f:=t^2+1\in {\mathcal W}\,$ since $\,V(f)\,$
contains $\,i\,$ and $\,j$; but $\,g:=(t-j)(t-i)\notin {\mathcal W}$, 
since $\,V(g)=\{i\}$.

\sk
(4) If $K$ is a field and $\,(S,D)=(I,0)$, the algebraic sets are precisely 
the {\it finite\/} subsets of $\,K$.  It follows that, in this case, the 
W-polynomials are precisely the completely split polynomials 
$\,(t-a_1)\cdots (t-a_n)\,$ where the $\,a_i\,$'s are distinct elements 
in $\,K$. 

\mk
(5) Let $\,(S,D)=(I,0)$, and let $\,a\in K\,$ be any element that is 
algebraic over $\,F=Z(K)$, with minimal polynomial $\,f(t)\in F[t]\,$ 
(in the usual field-theoretic sense).  By 
Wedderburn's Theorem in [We], the conjugacy class $\,\D:=\D(a)\,$ is 
algebraic (in our sense), with $\,f_{\D}=f(t)$. Therefore, $f\in F[t]\,$ 
is a W-polynomial over $\,K$. From this example, we see that the 
W-polynomials introduced in (3.1) are a generalization of the 
minimal polynomials of algebraic elements studied by Wedderburn in [We].

\bs
One reason the family $\,{\mathcal W}\,$ commands our interest is that
there are many possible ways to characterize the Wedderburn polynomials.
The characterizations in (3.2) above are only the simpler ones that follow
directly from the definitions. To close this section, we'll mention below 
a few other more substantial characterizations.

\bs
Consider for the moment the case where $\,K\,$ is a field and $\,(S,D)=(I,0)$.
For a monic polynomial  $\,f=t^n+b_{n-1}t^{n-1}+\cdots +b_0\in R$, we have
the following well-known ``companion matrix'':
$$  C(f):=\left( \begin{array}{ccccc}
0          &      1   & 0      &      \cdots        &    0   \\
0  &  0   &  1   &    \cdots          &  0    \\
\vdots   &  \vdots    &  \vdots    &   \ddots   &   \vdots     \\
0  &  0   &  0   &    \dots          &  1    \\
-b_0   &  -b_1   &   -b_2   &  \cdots   &  -b_{n-1}
\end{array}  \right),  \leqno (3.4)  $$
whose characteristic and minimal polynomials are both $\,f$.  For this
companion matrix to be {\it diagonalizable\/}, the classical criterion is 
that $\,f\,$ has $\,n\,$ distinct (characteristic) roots in $\,K$.
By (3.3)(4) above, this is {\it also\/} the criterion for $\,f\,$ to be a 
Wedderburn polynomial!  Through this example, we see that 
W-polynomials ought to be related to diagonalization questions in general.

\bs
In the case when $\,R=K[t,S,D]$, we can consider the $\,R$-module $\,R/Rf\,$
as a left $\,K$-vector space, with the $\,K$-basis $\,\{1,t,\dots,t^{n-1}\}$.
With respect to this basis, $\,C(f)\,$ is the matrix corresponding to 
left multiplication by $\,t$.  However, this action of $\,t\,$ is only a 
{\it pseudo-linear\/} transformation (in the sense that 
$\,t(kv)=S(k)t(v)+D(k)v$, for $\,k\in K\,$ and $\,v\in R/Rf$). Upon a change 
of $\,K$-basis on $\,R/Rf$, the matrix of a pseudo-linear transformation 
will change by an ``$(S,D)$-conjugation'' (analogous to the
$\,(S,D)$-conjugation of elements defined in \S2).  Using such a generalized 
notion of conjugation on matrices, we have the following characterization 
result for $\,W$-polynomials, which, in particular, extends the classical 
results on the diagonalizability of the companion matrix $\,C(f)$.

\bs\nt
{\bf Theorem 3.5.} {\it For any polynomial $\,f=t^n+b_{n-1}t^{n-1}+\cdots
+b_0\in R$, the following are equivalent\/}:

\sk
(1) {\it $f$ is a product of linear factors, and $\,R/Rf\,$ is a
semisimple $\,R$-module\/}.

\sk
(2) {\it The companion matrix $\,C(f)\,$ in $(3.4)$ can be 
$\,(S,D)$-conjugated to a diagonal matrix\/}.

\sk
(3) {\it There exists an invertible Vandermonde matrix\,\footnote{For any
given triple $\,(K,S,D)$, the general Vandermonde matrix
$\,V(c_1,\dots,c_n)\,$ is defined to be the $\,n\times n\,$ matrix whose 
$\,(i+1)\,$st row
is $\,\bigl(N_{i}(c_1),\dots, N_{i}(c_n)\bigr)$, where $\,N_i\,$ is the
$\,i\,$th power function (with respect to $\,(S,D)$) defined in (2.2).}
$\,V=V(c_1,\dots,c_n)\,$ over $\,K\,$ such that} 
$\,C(f)V=S(V)\,\mbox{diag}(c_1,\dots,c_n)+D(V)\,$.  {\it (Here, $S$ and
$D$ operate on matrices in the obvious way.)}

\sk
(4) deg$(f)=\sum_{j=1}^r\,\mbox{dim}\,_{C_j}\,E(f,a_j)$.

\bs
In this theorem, we have slipped in another characterization (4), which 
requires some notational explanations.  The elements $\,a_1,\dots,a_r\,$
there are representatives of the $\,(S,D)$-conjugacy classes of $\,K\,$ 
in which $\,f\,$ has roots, and the ``exponential space'' $\,E(f,a_j)\,$ is
defined to be 
$$\,\{0\}\cup \{x\in K^*:\;f(a_j^x)=0\}.$$  
This is a right vector space over the division subring
$$\,C_j:=\{0\}\cup \{c\in K^*:\;a_j^c=a_j\} \;\;\;\;\mbox{(the 
``$(S,D)$-centralizer'' of $\,a_j$)}.$$
For any $\,f$, it is true that
$\,\sum_{j=1}^r\,\mbox{dim}\,_{C_j}\,E(f,a_j)\leq \mbox{deg}(f)$; (4) of
the theorem above tells us that equality here is a criterion for $\,f\,$
to be a Wedderburn polynomial.

\sk

\bs
\begin{center}
{\bf \S4. \ Left-Right Symmetry of W-Polynomials}
\end{center}

\bs
According to the discussions in \S2,  W-polynomials are the (monic) 
llcm's of polynomials of the form $\,t-a\,$ ($a\in K$). In other words, 
a monic polynomial $\,f\,$ belongs to $\,{\mathcal W}\,$ iff the left ideal
$\,Rf\,$ can be expressed as an intersection of the form $\,\bigcap \,R(t-a)$,
where $\,a\,$ ranges over some subset of $\,K$. From this description,
it would appear that $\,{\mathcal W}\,$ is a family dependent on the
left structure of the ring $\,R$.  Our first result on $\,{\mathcal W}\,$
is that it is in fact a left-right symmetric family.

\bs
To state this symmetry result more formally, let us introduce the following
notations.  Let $\,{\mathcal W}^{\ell}\,$ be the family of monic polynomials 
$\,f\in R\,$ such that $\,Rf\,$ is an (arbitrary) intersection of principal 
left ideals of the form $\,R(t-a)\,$ ($a\in K$), and let 
$\,{\mathcal W}^{\ell}_0\,$ be the family of monic polynomials $\,f\in R\,$ 
such that $\,Rf=R(t-a_1)\cap \cdots \cap R(t-a_n)$, where $\,a_i\in K\,$ and 
$\,n=\mbox{deg}(f)$.  By using right (instead of left) principal ideals, 
we can similarly define two families of monic polynomials $\,{\mathcal W}^r\,$
and $\,{\mathcal W}^r_0\,$. The symmetry result [LL$_4$:~(4.5)] states the 
following:

\bs\nt
{\bf Theorem 4.1.} 
$\,{\mathcal W}={\mathcal W}^{\ell}={\mathcal W}^{\ell}_0
={\mathcal W}^{r}_0={\mathcal W}^{r}$.

\bs
The first two equalities here are already covered by our previous discussions,
so the gist of the theorem is in the last two equalities.  We refer the reader 
to [LL$_4$] for the proof of (4.1), but let us give here some indication of 
how one could achieve a passage from the left structure of $\,R\,$ to its 
right structure. The following result, in fact, gives a {\it constructive\/} 
description of the inclusion 
$\,{\mathcal W}^{\ell}_0 \subseteq {\mathcal W}^{r}_0\,$ in (4.1).

\bs\nt
{\bf Theorem 4.2.} {\it Suppose $\,f\in R\,$ is monic of degree $\,n$, 
and $\,Rf\,$ has a representation as $\,R(t-a_1)\cap \cdots \cap R(t-a_n)$.}
({\it In other words, $\,f\,$ is the minimal polynomial of the P-independent
set $\,\{a_1,\dots,a_n\}$.}) {\it If $\,h_i\,$ is the minimal polynomial of 
$\,\{a_1,\dots,a_{i-1},a_{i+1},\dots,a_n\}$, then, for}
$\,b_i:=a_i^{h_i(a_i)}\;(1\leq i\leq n)$, {\it we have a
representation $\,fR=\bigcap_{i=1}^n (t-b_i)R$.}

\bs
A constructive description of the reverse inclusion
$\,{\mathcal W}^{r}_0 \subseteq {\mathcal W}^{\ell}_0\,$
can be given ``similarly''.

\bs
The Symmetry Theorem (4.1) serves to show that the family $\,{\mathcal W}\,$ 
of Wedderburn polynomials is sufficiently intrinsic to the ring $\,R\,$
to be characterizable in terms of either its left structure or its right 
structure. In the case when $\,S\,$ is an automorphism of $\,K$, this is
not too surprising, and was in fact known in a more general form to Ore 
(see [Or:~Ch.~II, Th.~4]).  But if $\,S\,$ is just an endomorphism of
$\,K$, the situation is different.  In case $\,S(K)\neq K$, $\,R\,$ 
is {\it not\/} a principal right ideal domain (or even a right Ore domain), 
so there is an apparent disparity between the left and the right structures 
of $\,R$. Here, one seems to have no {\it \`a priori\/} reason to expect a 
symmetry result such as (4.1).

\bs
What comes to one's rescue is the notion of a {\it semifir\/} due to 
P.~M.~Cohn. A semifir is a domain in which every finitely generated left ideal 
is free of a unique rank --- and this is known to be a left-right symmetric 
notion (see [Co$_2$]).  Since our ring $\,R=K[t,S,D]\,$ is a principal left 
ideal domain, it is a (left and and hence right) semifir.  From this fact,
it turns out that one has enough symmetry information on $\,R\,$ to prove 
Theorem 4.1. We refer the reader to [LL$_4$] for the details of the proof.

\bs
An interesting byproduct of these considerations is a theory of {\it left\/} 
roots of polynomials in $\,R$.  For $\,f\in R$, we have defined earlier
$$ V(f)=\{a\in K:\;\,f(a)=0\} = \{a\in K:\;\,f\in R(t-a)\};  \leqno (4.3) $$
such elements $\,a\,$ should be called the {\it right roots\/} of $\,f$.
We can correspondingly define 
$$ V'(f)=\{b\in K:\;\,f\in (t-b)R\}; \leqno (4.4) $$
such $\,b\,$'s may then be called the {\it left roots\/} of $\,f$.  This 
opens the way to many more interesting facts.  For instance, in analogy to
(3.2)(3), one can prove the following characterization of W-polynomial
in terms of left root sets:

\bs\nt
{\bf Proposition 4.5.} {\it For any monic $\,f$, we have
$\,f\in {\mathcal W}\,$ iff, $\,\forall\;q\in R$, $\,V'(f)\subseteq V'(q)$
\linebreak
implies $\,q\in fR$.}

\bs 
In general, however, the right and left root sets $\,V(f)\,$ and $\,V'(f)\,$
behave rather differently, especially in the case when $\,S(K)\neq K$.
In this case, the theorem below says that $\,V'(f)\,$ is situated in a 
very special way with respect to the coset partition of $\,K\,$ modulo the 
additive subgroup $\,S(K)$.  This result has no right root analogue,
and thus represents a new feature of the left root theory.

\bs\nt
{\bf Theorem 4.6.} {\it For any nonzero polynomial $\,f$, either all left 
roots of $\,f\,$ lie in a single $\,S(K)$-coset of $\,K$, or no two left roots
lie in the same $\,S(K)$-coset. If $\,f\,$ is monic, then necessarily all 
left roots of $\,f\,$ lie in a single $\,S(K)$-coset of $\,K$.}

\sk

\bs
\begin{center}
{\bf \S5. \ Factor Theorem, Product Theorem, and the $\,\Phi$-Transform}
\end{center}

\bs
Wedderburn polynomials turn out to have very tractable properties under
multiplication and division.  In this section, we'll state several theorems 
on the factors and products of W-polynomials.

\bs
First we have to clarify what exactly is meant by the word ``factor'' in this 
paper.  Throughout the sequel, we'll say that a polynomial $\,p\,$ is a 
{\bf factor} of another polynomial $\,f\,$ if $\,f= p_1pp_2\,$ for some
$\,p_1,\,p_2 \in R$.  Right and left 
factors of $\,f\,$ have their usual meanings, and these are, of course, 
particular kinds of factors in our sense.  The main result about factors
of W-polynomials is the following (cf.~[Or:~Ch.\,II, Th.\,3]):

\bs\nt
{\bf Factor Theorem 5.1.} {\it For any monic $\,f\in R$, the following are 
equivalent\/}:

\sk
(1) $\,f\,$ {\it is a W-polynomial}\;;

\sk
(2) $\,f\,$ {\it splits completely, and every monic factor of $\,f\,$ is 
a W-polynomial}\;;

\sk
(3) $\,f\,$ {\it splits completely, and every monic quadratic factor of 
$\,f\,$ is a W-polynomial\/}.

\bs
Next, we'll discuss the situation of {\it products.} In the classical case 
when $\,(S,D)=(I,0)\,$ and $\,K\,$ is a field, we see already that the product
of two W-polynomials need not be a W-polynomial.  In fact, if $\,g,\,h\in
{\mathcal W}\,$ in this case, we'll have $\,gh\in {\mathcal W}\,$ iff
$\,V(g)\cap V(h)=\emptyset$.  One must then look for generalizations
of this statement to the case of $\,R=K[t,S,D]\,$ (for any triple
$\,(K,S,D)$).  

\bs
To solve this problem, we need the tool of a certain ``$\Phi$-transform''.
In the Product Formula (2.3) for evaluating $\,gh\,$ at $\,a$, we first 
encountered the expression $\,a^{h(a)}\,$ (in case $\,h(a)\neq 0$);
similar expressions also occurred in the statement of Theorem (4.2).
This led us to the following formal definition.

\bs\nt
{\bf Definition 5.2.} For $\,h\in R=K[t,S,D]$, we define the 
``$\Phi\,$-transform'' (associated to $\,h$)
$$\,\Phi_h : \;K\setminus V(h) \longrightarrow K   \leqno (5.3)$$ 
by $\,\Phi_h(x)=x^{h(x)}$, whenever $\,h(x)\neq 0$.  (We do not attempt to 
define $\,\Phi_h\,$ on $\,V(h)$.)  Note that $\,\Phi_h\,$ always preserves 
the $\,(S,D)$-conjugacy class of its argument $\,x$.  

\bs\nt
{\bf Examples 5.4.}  

\mk
(1) Say $\,h\,$ is a nonzero constant (polynomial) 
$\,c\in K^*$.  Then, $\,\Phi_c(x)=x^c\,$ for all $\,x\in K$, so 
$\,\Phi_c\,$ is defined on all of $\,K\,$ and is exactly $\,(S,D)$-conjugation
by the element $\,c$. (In particular, $\,\Phi_1\,$ is just the identity map 
on $\,K$.)  In view of this example, we can think of the $\,\Phi$-transform
as a kind of generalization of $\,(S,D)$-conjugation.

\sk
(2) In the case $\,D=0\,$ and $\,h=t$, $\,\Phi_t(a)=a^a=S(a)aa^{-1}=S(a)\,$
for any $\,a\neq 0$.  Here, $\,\Phi_t\,$ is the restriction of $\,S\,$
to $\,K^*$.

\bs
With the aid of the $\,\Phi$-transform, one can prove the following
result which gives various criteria for the product of two W-polynomials
to be another W-polynomial.

\bs\nt
{\bf Product Theorem 5.5.} {\it For $\,f=gh\in R\,$ where $\,g,\,h\,$ are
monic, the following are equivalent\/}:

\mk
(1) $\,f\in {\mathcal W}$;

\sk
(2) $\,g,\,h\in {\mathcal W}$, {\it and} $\,1\in Rg+hR$;

\sk
(3) $\,g,\,h\in {\mathcal W}$, {\it and} $\,V(g)\subseteq \mbox{im}(\Phi_h)$;

\sk
(4)  $\,g,\,h\in {\mathcal W}$, {\it and, $\forall\,a\in V(g),\;b\in V'(h)$,
we have $\,(t-a)(t-b)\in {\mathcal W}$}.

\bs
The formulation (4) of a W-polynomial criterion in terms of {\it quadratics\/}
is particularly nice since (3.3)(3) makes it quite easy to recognize quadratic
W-polynomials; see also Theorem (8.2) below.

\sk

\bs
\begin{center}
{\bf \S6. \ Three Rank Theorems}
\end{center}

\bs
In this section, we state a few theorems concerning the ranks of algebraic 
sets.  Recall that, for an algebraic set $\,\Delta$, rk$\,(\Delta)\,$ is the 
degree of its minimal polynomial $\,f_{\D}$, and the {\it P-closure\/} 
$\,\overline{\D}\,$ of $\,\D\,$ is (or may be taken to be) $\,V(f_{\D})$. 
We'll say that an algebraic set $\,\D\,$ is {\it full\/} if 
$\,\overline{\D}=\D$.

\bs
The first rank formula below gives a relation between the ranks of the
union and intersection of algebraic sets.  The formula is obviously
prompted by the well-known dimension equation in the theory of 
finite-dimensional vector spaces. However, the use of the P-closures
in the intersection $\,\overline{\D}\cap \overline{\G}\,$ below adds a
subtle element to the formula. (Of course, the formula would look simpler
if we restrict ourselves to {\it full\/} algebraic sets.  But this would
be an unnecessary sacrifice of generality.)

\bs\nt
{\bf First Rank Theorem 6.1.} {\it For any two algebraic sets $\,\D\,$ and 
$\,\G$, we have}
$$ \mbox{rk}\,(\D)+\mbox{rk}\,(\G) = \mbox{rk}\,(\D\cup \G) + 
\mbox{rk}\,\bigl(\overline{\D} \cap \overline{\G} \bigr).  $$
{\it In particular,} $\,\mbox{rk}\,(\D\cup \G)
=\mbox{rk}\,(\D)+\mbox{rk}\,(\G)\,$ {\it iff\/} 
$\,\overline{\D}\cap \overline{\G}=\emptyset $.

\bs
We move on now to the Second Rank Theorem, which deals with the change 
of the rank of an algebraic set under a $\,\Phi$-transform, say $\,\Phi_h$.
One basic example to keep in mind is when $\,h\,$ is a constant polynomial
$\,c\in K^*$.  In this case, $\,\Phi_h\,$ is $\,(S,D)$-conjugation by
$\,c\,$ (see (5.4)(1)), which is easily seen to be {\it rank-preserving\/} 
on algebraic sets.  In the general case, a $\,\Phi$-transform $\,\Phi_h\,$ 
will be {\it rank-decreasing\/}, and the precise result is as follows.

\bs\nt
{\bf Second Rank Theorem 6.2.} {\it Let $\,h\in R$, and $\,\D\subseteq K\,$ 
be an algebraic set disjoint from $\,V(h)$. Then}
$$\,\mbox{rk}\,\bigl(\Phi_h(\D)\bigr) = \mbox{rk}\,(\D) - 
\mbox{rk}\,\bigl(\overline{\D}\cap V(h) \bigr).$$
{\it In particular,} $\,\mbox{rk}\,\bigl(\Phi_h(\D)\bigr)=\mbox{rk}\,(\D)\,$
{\it iff\/}  $\,\overline{\D}\cap V(h)=\emptyset$.

\bs
A powerful application of the rank formula for $\,\Phi_h\,$ is the 
following result giving a natural bound on the rank of the zero set of a 
product of two arbitrary polynomials.

\bs\nt
{\bf Third Rank Theorem 6.3.} {\it For $\,f=gh\in R\setminus \{0\}$,
we have} 
$$\,\mbox{rk}\,(V(f))\leq \mbox{rk}\,(V(g))+\mbox{rk}\,(V(h)).$$

\mk
In spite of its simplicity, (6.3) is not an easy result. Note, for
instance, the following interesting application of it to W-polynomials. 
Suppose $\,g,h\,$ above are monic, of degrees $r$ and $s$ respectively. 
If we are given $\,f:=gh\in {\mathcal W}$, then by (3.2), 
$\,\mbox{rk}\,(V(f))=\mbox{deg}(f)=r+s$, and (6.3) gives 
$$\,\mbox{rk}\,(V(g))+\mbox{rk}\,(V(h))\geq r+s=
\mbox{deg}(g)+\mbox{deg}(h).$$
This implies that $\,\mbox{rk}\,(V(g))=r\,$ and $\,\mbox{rk}\,(V(h))=s$,
and hence (again by (3.2)) $\,g,\,h\in {\mathcal W}$.  This conclusion
amounts essentially to the implication $(1)\Longrightarrow (2)$ in
the Factor Theorem (5.1).

\sk

\bs
\begin{center}
{\bf \S7. \ A Duality Between Two Complete Modular Lattices}
\end{center}

\bs
The considerations in the previous sections lead us quickly to 
the construction of several lattices.  There is only one more piece 
missing in the whole puzzle.  This is given by the following proposition,
where ``rgcd'' is used as an abbreviation for ``right greatest
common divisor''.

\bs\nt
{\bf Proposition 7.1.} {\it The intersection of any nonempty family of full 
algebraic sets $\,\{\D_j:\,j\in J\}\,$ is also a full algebraic set, with 
minimal polynomial given by} rgcd$\,\{f_{\D_j}:\,j\in J\}$.

\bs\nt
{\bf Remark 7.2.} The fact that $\,\bigcap_j\,\D_j\,$ has minimal polynomial
$\,\mbox{rgcd}\{f_{\D_j}:\,j\in J\}\,$ is generally not true if the algebraic
sets $\,\D_j\,$ are not all {\it full\/}.  For instance, in the division ring 
$\,K\,$ of real quaternions, $\,\D=\{i\}\,$ is full and $\,\G=\{j,\,k\}\,$ is 
not full.  The rgcd of $\,f_{\D}=t-i\,$ and $\,f_{\G}=t^2+1\,$ is $\,t-i$.
But $\,\D\cap \G=\emptyset\,$ has minimal polynomial $\,1$.

\bs
Now consider the poset $\,{\mathcal F}={\mathcal F}(K,S,D)\,$ of all 
{\it full\/} algebraic sets in $\,K\,$ (with respect to $\,(S,D)$), where 
the partial ordering is given by inclusion:
$$\,\D\leq \G \;\Longleftrightarrow \; \D\subseteq \G  \;\;\,
(\mbox{for} \,\; \D,\,\G \in {\mathcal F}\,). \leqno (7.3) $$ 
This poset is a {\it lattice\/}, with $\,\D \wedge \G\,$
given by $\,\D\cap \G\,$ (which lies in $\,{\mathcal F}\,$ by (7.1)), and
with  $\,\D \vee \G \,$ given by $\,\overline{\D\cup \G}$.  (The union
$\,\D\cup \G\,$ is algebraic, but may not be {\it full\/}, as the case
$\,|\D|=|\G|=1\,$ already shows.)  Note that the lattice $\,{\mathcal F}\,$ 
has a smallest element, given by the empty set $\,\emptyset $ 
(see (3.3)(1)).

\bs
If $\,K\,$ itself is not $\,(S,D)$-algebraic, $\,{\mathcal F}\,$ will not 
have a largest element.  In this case, we can define an augmented lattice 
$\,{\mathcal F}^*\,$ by adjoining to $\,{\mathcal F}\,$ a largest element.
If $\,K\,$ is $\,(S,D)$-algebraic, then it is the largest element in
$\,{\mathcal F}$, and we can simply define $\,{\mathcal F}^*\,$ to be 
$\,{\mathcal F}$.

\bs
To get a lattice in the context of (skew) polynomials, we consider the set 
$\,{\mathcal W}={\mathcal W}(K,S,D)\,$ of all 
W-polynomials, partially ordered by: 
$$  f\leq h \; \Longleftrightarrow \; Rf\subseteq Rh \;
          \Longleftrightarrow \; h\;\,\mbox{is\;a\;right\;divisor\;of}\;\,f
\;\;\;(\mbox{for}\;\,f,\,h\in {\mathcal W}).  \leqno (7.4)  $$
The poset $\,{\mathcal W}\,$ is again a lattice: for 
$\,f,\,h\in {\mathcal W}\,$ as above, $\,f \vee h \,$ is given
by $\,\mbox{rgcd}(f,\,h)\,$ (this being a W-polynomial by the Factor
Theorem), and $\,f \wedge h\,$ is given by
$\,\mbox{llcm}(f,\,h)\,$ (this being a W-polynomial since it is the minimal
polynomial of $\,V(f)\cup V(h)$).  The lattice $\,{\mathcal W}\,$
has a largest element, given by $\,1\in {\mathcal F}$, and it will have 
a smallest element iff $\,K\,$ happens to be $\,(S,D)$-algebraic. (The 
smallest element in the latter case is the minimal polynomial of $\,K\,$ 
itself.) In analogy with the case of full algebraic sets, we can introduce 
an {\it augmented W-polynomial lattice\/} $\,{\mathcal W}^*$, which is 
defined to be $\,{\mathcal W}\,$ if $\,K\,$ is $\,(S,D)$-algebraic, and 
$\,{\mathcal W}\,$ adjoined with the polynomial $\,0\,$ otherwise (this 
being smaller than all other W-polynomials).

\bs
Once we have set up the augmented lattices $\,{\mathcal F}^*\,$ and
$\,{\mathcal W}^*\,$ as above, we can use (6.1), (7.1), etc. to prove the 
following duality result.

\bs\nt
{\bf Theorem 7.5.} {\it For a fixed triple $\,(K,S,D)$, we have the 
following\/}:

\sk
(1) {\it ${\mathcal F}^*\,$ and $\,{\mathcal W}^*\,$ are both complete 
modular lattices.  The maps $\,\D\mapsto f_{\D}\,$ and $\,f\mapsto V(f)\,$ 
(extended in the obvious way) define mutually inverse lattice dualities
between} $\,{\mathcal F}^*\,$ and $\,{\mathcal W}^*$.  

\sk
(2) {\it The map\/} ``rk'' {\it (extended in the obvious way) is the 
``dimension function'' on the modular lattice $\,{\mathcal F}^*\,$ (in the 
sense of lattice theory), and the degree map\/} ``deg'' {\it is the ``dual 
dimension function'' on the modular lattice $\,{\mathcal W}^*$.}

\sk
(3) {\it The minimal elements (the so-called atoms) of the lattice 
$\,{\mathcal F}^*\,$ are the singleton subsets of $\,K$, and the maximal 
elements of the lattice $\,{\mathcal W}\,$ are the monic linear polynomials 
in $\,R$.}

\sk
(4) {\it For $\,f\leq h\in {\mathcal W}$, the ``interval'' $\,[f,\,h]\,$
in the lattice $\,{\mathcal W}\,$ is isomorphic to the lattice of
all $\,R$-submodules of $\,Rh/Rf$.}

\bs
In our proof of (1) of this theorem in [LL$_5$], the modularity of 
$\,{\mathcal F}^*\,$ is {\it not\/} proved directly, but is rather deduced 
from the modularity of $\,{\mathcal W}^*$.  It is, therefore, of interest 
to record the following statement, which essentially amounts to the modular 
law for $\,{\mathcal F}^*$.

\bs\nt
{\bf Proposition 7.6.} {\it Let $\,\G,\;\Pi\,$ and $\,\D\,$ be algebraic sets, 
where $\,\G,\;\Pi\,$ are full, and $\,\D\subseteq \G$.  If $\,x\in \G\,$ 
is P-dependent on $\,\Pi\cup \D\,$, then it is already P-dependent on the 
smaller set $\,(\G\cap \Pi)\cup \D$.}

\sk

\bs
\begin{center}
{\bf \S8. \ The $\,(S,D)$-Metro Equation}
\end{center}

\bs
In the theory of division rings, the study of the equation $\,ax-xb=c\,$
has had a long history, going back to the work of R.~E.~Johnson [Jo] and
N.~Jacobson [Ja$_2$].  By a slight abuse of the terminology of
P.~M.~Cohn ([Co$_1$], [Co$_3$]), we shall call $\,ax-xb=c\,$ a 
``metro equation'' over $\,K$.  It turns out that the notion of metro
equations bears a close relationship to that of Wedderburn polynomials.

\bs
First, we generalize the classical metro equation notion to our 
$\,(S,D)$-setting. This is not difficult: for $\,a,b,c\in K$, let us call
$$ax-S(x)b-D(x)=c  \leqno (8.1)$$
the $\,(S,D)$-{\it metro equation\/} (associated with $\,a,b,c$).  (Of course,
when $\,(S,D)=(I,0)$, (8.1) boils down to the ordinary metro equation
$\,ax-xb=c$.) In the case $\,c=0$, (8.1) has an obvious solution $\,x=0$, 
so in the following, we'll assume $\,c\neq 0\,$ whenever (8.1) is considered.
The following result gives the precise relationship between (8.1) and 
{\it quadratic\/} Wedderburn polynomials, in the general $\,(S,D)$-setting.

\bs\nt
{\bf Theorem 8.2.} {\it For any $\,a,\,b\in K\,$ and $\,c\in K^*$, the
following are equivalent\/}:

\sk
(1) {\it The $\,(S,D)$-metro equation $\,ax-S(x)b-D(x)=c\,$ has a solution
in $\,K\,$};  

\sk
(2) $(t-b^c)(t-a)\in {\mathcal W}\,$;

\sk
(3) $1\in R\cdot (t-b^c)+(t-a)\cdot R.$

\bs
Here, $(2)\Longleftrightarrow (3)$ follows directly from (5.5), so what is 
new in (8.2) is the equivalence $(1)\Longleftrightarrow (2)$. Just for the 
sake of illustration, let us give the proof for this equivalence in the
``classical case'', namely, when $\,(S,D)=(I,0)$. This will, in
fact, serve as the {\it one and only\/} proof given in this entire paper.

\bs
What we are trying to prove here is that $\,ax-xb=c\,$ has a solution in 
$\,K\,$ iff the quadratic
$$\,f(t):=\bigl(t-b^c\bigr)(t-a)=t^2-\bigl(a+b^c\bigr)t+b^ca \in K[t] $$
(where now $\,b^c:=cbc^{-1}$) is a W-polynomial.  For the latter to be 
true, we need to have 
a (right) root of $\,f\,$ in $\,K\,$ that is different from $\,a\,$ (see 
(3.3)(3)).  Write such a root (if it exists) in the form $\,a-cx^{-1}\,$ 
(where $\,x\in K^*$).  For this to be a right root of $\,f$, we need:
\begin{eqnarray*}
0=f\bigl(a-cx^{-1}\bigr)&=&(a-cx^{-1})^2-\bigl(a+b^c\bigr)(a-cx^{-1})+b^ca \\
&=& \bigl(a-cx^{-1}-a-b^c\bigr)(a-cx^{-1})+b^ca \\
&=& \bigl(-cx^{-1}-b^c\bigr)(a-cx^{-1})+b^ca \\
&=& c\bigl(-x^{-1}a+bx^{-1}+x^{-1}cx^{-1}\bigr),
\end{eqnarray*}
that is, $\,x^{-1}a-bx^{-1}=x^{-1}cx^{-1}$. Multiplying this equation by 
$\,x\,$ from the left and from the right, we get it down to $\,ax-xb=c$. 
Therefore, finding a root for $\,f\,$ other than $\,a\,$ amounts exactly 
to solving the metro equation $\,ax-xb=c\,$!

\bs\nt
{\bf Example 8.3.} Let $\,K\,$ be the division hull of the Weyl algebra
$\,{\mathbb C}\langle u,v\rangle\,$, with the relation $\,uv-vu=1$.
Choosing $\,a=b=u$, we see that the metro equation $\,ax-xb=1\,$ has
a solution $\,x=v$.  Thus, it follows from (8.2) that $\,f(t):=(t-u)^2\,$
is a W-polynomial over $\,K$.  (Two roots for $\,f\,$ are $\,u\,$ and
$\,u-v^{-1}$, according to the calculation above.) This example is 
interesting, as it shows that a Wedderburn polynomial may very well have 
{\it repeated\/} linear factors.

\bs
Theorem (8.2) has some nice applications to the solvability of the metro 
equation (8.1) in case the $\,(S,D)$-conjugacy class of $\,b\,$ is an 
algebraic set.  We close by stating the following consequence of (8.2).

\bs\nt
{\bf Proposition 8.4.} {\it Suppose the $\,(S,D)$-conjugacy class 
$\,\D:=\D^{S,D}(b)\,$ is algebraic, and let $\,a\in K\setminus \D$, 
$\,c\in K^*$. Then $\,(t-b^c)(t-a)\in {\mathcal W}$, and the $\,(S,D)$-metro
equation\/} (8.1) {\it has a unique solution.}

\sk

\bs

\end{document}